\documentclass[11pt]{amsart}

\usepackage[a4paper,margin=1in]{geometry}
\usepackage{amsmath,amsrefs, amssymb,amsthm,mathtools}
\usepackage[pdfencoding=auto,psdextra,colorlinks=true,linkcolor=black]{hyperref}
\usepackage{graphicx}
\usepackage{amscd,amsfonts,graphicx,latexsym,enumitem,imakeidx}
\usepackage{tikz-cd}
\usepackage{blindtext}
\usepackage{listings}
\usepackage[all,arc]{xy}
\usepackage{mathrsfs}
\usepackage{cite}
\usepackage{xcolor}

\title{Elementary Factorization for Continuous Orthogonal Matrices}

\author[G. Huang]{Gaofeng Huang}
\address{Department of Mathematics \\ Harvard University \\  
Cambridge, MA 02138, USA}
\email{ghuang@math.harvard.edu}

\author[B. Tran]{Phan Quoc Bao Tran}
\address{Mathematisches Institut \\ University of Bern \\  
Bern, Switzerland}
\email{phan.tran@unibe.ch}

\newtheorem{theorem}{Theorem}[section]
\newtheorem{lemma}[theorem]{Lemma}
\newtheorem{proposition}[theorem]{Proposition}
\newtheorem{corollary}[theorem]{Corollary}
\theoremstyle{definition}

\theoremstyle{remark}
\theoremstyle{remark}
\newtheorem{remark}[theorem]{Remark}

\newcommand{\C}{\mathbb{C}}

\newcommand{\diag}{\operatorname{diag}}


\thanks{The first author was funded by Schweizerischer Nationalfonds   Postdoc.Mobility [P500-$2\_239113$]. The second author was funded by Schweizerischer Nationalfonds  [10005432] and [229554].}

\begin{document}
\maketitle
\begin{abstract}
    We prove that every nullhomotopic continuous map, from a finite-dimensional
    normal topological space $X$ to the complex special orthogonal group,  is a product of elementary orthogonal
    matrices with entries in the ring of continuous functions. The
    number of factors is bounded by a constant depending only on the size of the matrices and the
    covering dimension of $X$. The same factorization result holds for nullhomotopic special orthogonal matrices with entries in a unital commutative Banach algebra.
\end{abstract}

\section{Introduction}

Let $R$ be a commutative ring with unit. A recurring theme originating in algebraic $K$-theory is to determine when a determinant-one matrix with entries in $R$ can be written as a finite product of elementary matrices $I_n+a E_{ij}$, where $i\ne j$ and $a\in R$, or equivalently as a finite product of unipotent lower and upper triangular matrices in alternating order. We call these unipotent triangular factors \emph{unitriangular} matrices. 

\smallskip
For rings of functions this becomes a factorization problem depending on parameters. Let $\mathcal{C}(X)$ be the ring of complex-valued continuous functions on a finite-dimensional normal topological space $X$. Leonid Vaserstein \cite{04073958} proved that a matrix in $SL_n(\mathcal{C}(X))$ admits such a factorization if and only if, viewed as a map $X\to SL_n(\C)$, it is homotopic to the constant identity map. We call such a map \emph{nullhomotopic}. 
Using Vaserstein's continuous factorization as topological input, Bj{\" o}rn Ivarsson and Frank Kutzschebauch \cite{06014072} solved the corresponding problem posed by Mikhail Gromov \cite{MR1001851} (for $X = \C^N$) for the ring $\mathcal{O}(X)$ of holomorphic functions on a reduced Stein space $X$. Namely every nullhomotopic holomorphic map $f: X \to SL_n(\C)$ is a product
\begin{align*}
    f(x) = \left(\begin{matrix} 1 & 0 \cr G_1(x) & 1 \cr \end{matrix} \right)   
\left(\begin{matrix} 1 & G_2(x) \cr 0 & 1 \cr \end{matrix} \right)  \ldots \left(\begin{matrix} 1 & G_L(x)\cr 0 & 1 \cr \end{matrix} \right) 
\end{align*}
where the $G_i$ are maps $G_i : X \to \mathbb{C}^{n(n-1)/2}$. A map admitting such a decomposition must be nullhomotopic: we simply rescale $G_i(x)$ by $t G_i(x), 0 \le t \le 1$. 

The main theoretical tool Ivarsson and Kutzschebauch used is the Oka principle for sections of stratified elliptic submersions  developed by Franc Forstneri{\v c} and Jasna Prezelj \cite[\S 7]{MR1812449}. A crucial step in \cite{06014072} is to reduce $n$ equations, describing the fiber of a certain map $\Phi$, to a single equation.  

In dimension one, Manfred Klein and Karl Ramspott \cite[\S IV]{MR0966022} showed that every holomorphic map from a noncompact Riemann surface to $SL_n(\C)$ is a product of unitriangular factors. 

\smallskip
For the symplectic group $Sp_{2n}(\C)$, the continuous case was settled by Ivarsson, Kutzschebauch and Erik L{\o}w \cite{MR4078081}, the holomorphic case with $n=2$ by the same authors \cite{MR4578529} and general $n$ by Josua Schott \cite{900569216}, there the $2n$ equations can only be reduced to $n$ equations which makes the argument much more technical. 
A refinement of Schott's result, allowing a new type of unitriangular factors, was obtained by Huang, Kutzschebauch and Tran in \cite{MR5101717}. With more freedom in choosing the factors, a reduction to one equation became possible. This motivates our choice of elementary factors in Section \ref{sec-eleMat}.

We refer to the survey  \cite{MR5024280} by Kutzschebauch and the first author for a broader picture about parametric matrix factorization.

\medskip
The purpose of this paper is to carry out the continuous part of this program for the orthogonal group, following Ivarsson, Kutzschebauch and L{\o}w \cite{MR4078081}. Fix an integer $m\ge 3$ and set
\begin{equation*}
S_{2n} = \begin{pmatrix}0&I_n\\ I_n&0\end{pmatrix}\ \text{if } m=2n,
\quad 
S_{2n+1}=\begin{pmatrix}0&0&I_n\\ 0&1&0\\ I_n&0&0\end{pmatrix}\ \text{if } m=2n+1,
\end{equation*}
and for a commutative ring $R$ with unit write
\begin{center}
$    SO_m(R,S_m):=\{ M \in SL_m(R): M^{T}S_m M=S_m \}.$
\end{center}
We will specify the elementary orthogonal factors in Section \ref{sec-eleMat}. 
Our main result is the following. 
\begin{theorem} \label{cont-factorization}
    For any natural number $m \ge 3$ and any integer $d \ge 0$ there exists a natural
    number $L = L(m, d)$ such that for any normal topological space $X$ of finite covering dimension $d$, every nullhomotopic continuous map $f : X \to SO_{m}(\mathbb{C},S_m)$ can be written as a product of at most $L$ elementary orthogonal matrices with entries in $\mathcal{C}(X)$. 
\end{theorem}

Here, the condition $m \geq 3$ is essential. For $m=2$, $SO_2(\mathbb C,S_2)$ is isomorphic to $\C^*$ and the only elementary orthogonal matrix is $I_2$. There is no elementary factorization for $SO_2(\mathbb C,S_2)$. 

Over an arbitrary field $K$, the so-called spinor norm is an obstruction to the factorization of the diagonal matrix 
$$           \diag(1,\dots,1,a,\,1,\dots,1,a^{-1},\cdots,1)
$$ 
if $a$ is not a square in $K$ \cite{MR0360645}. Over the complex numbers, $a$ is a square and factorization is always possible. We prove in Section \ref{subsect-Whitehead} that such a diagonal matrix with $a$ a square can be factorized in $SO_m(\mathcal{C}(X),S_m)$.  

\subsection{Proof strategy of Theorem \ref{cont-factorization}} 

Following \cite{MR4078081} the proof has three main ingredients. First, we use the Gram--Schmidt orthogonalization process (Proposition \ref{theo-GS}) to change $f$ to a nullhomotopic continuous map $X \to SO_{m}(\mathbb{C}, S_m) \cap U(m)$, which is a compact subgroup (Proposition \ref{K_U-compact}). Thanks to this compactness, we can apply the Calder--Siegel theorem (Proposition \ref{Cardel-Siegel}) to get a uniform homotopy $H_t: X \to SO_{m}(\mathbb{C}, S_m) \cap U(m)$. That is, $H_1$ is a finite product of factors each uniformly close to $I_{m}$. On the other hand, the Gauss--Jordan elimination process (Proposition \ref{lemm-GJ-orth}) for continuous orthogonal matrices near the identity gives a factorization into elementary orthogonal factors. Combining these decompositions together, we get the desired unitriangular factorization for the nullhomotopic continuous map $f : X \to SO_{m}(\mathbb{C}, S_m)$.


As a byproduct of the Gauss--Jordan process, Proposition \ref{lemm-GJ-orth}, we get a factorization result for unital commutative Banach algebra $R$ similar to \cite[Theorem 1.1]{MR4078081}, cf.\@ Milnor \cite[Lemma 7.1]{MR0349811}. 

\begin{theorem}
For a unital commutative Banach algebra $R$, let $M\in SO_{m}(R,S_m)$ be nullhomotopic
in $SO_{m}(R,S_m), m \ge 3$. Then $M$ can be written as a product of elementary orthogonal matrices in $SO_{m}(R,S_m)$. 
\end{theorem}

The paper is organized as follows. We define the elementary orthogonal matrices in Section \ref{sec-eleMat} and explain the factorization for the diagonal type matrix. Next in Section \ref{sec-GJ} we show the Gauss--Jordan elimination for matrices near the identity over a unital commutative  Banach algebra. As an application we have in Section \ref{sec-Banach} the factorization of special orthogonal matrices over unital commutative Banach algebras. We perform in Section \ref{sec-GS} the Gram--Schmidt process and apply a result by Calder--Siegel to obtain a uniform nullhomotopy. Section \ref{sec:combine} threads the pieces together into a proof of Theorem \ref{cont-factorization}.

\section{Elementary orthogonal matrices} \label{sec-eleMat}
\subsection{Even case}
Let $E_{ij}$ be the elementary $n\times n$ matrix with $1$ in position $(i,j)$ and zeros elsewhere.
For a matrix $M\in M_{2n}(R)$ denote the rows of $M$ by $v_1,\dots,v_n,w_1,\dots,w_n$.

We introduce the following \emph{even elementary orthogonal matrices} in $SO_{2n}(R, S_{2n})$.

For $i\neq j$ and $a \in R$, set 
\begin{align*}
    K_{ij}(a):= \begin{pmatrix}
                    I_n+aE_{ij} & 0\\
                    0 & I_n-aE_{ji}
                \end{pmatrix}
\end{align*} 
Multiplying $M$ by $K_{ij}(a)$ from the left replaces the $i$th row $v_i$ by $v_i+a v_j$ and the $(n+j)$th row $w_j$ by $w_j-a w_i$.

\smallskip

For $1\le i \neq j\le n$ and $a\in R$, set
\begin{center}
$F_{ij}(a):=
\begin{pmatrix}
I_n & 0\\
a(E_{ij}-E_{ji}) & I_n
\end{pmatrix}, $
\end{center}
Multiplying $M$ by $F_{ij}(a)$ from the left replaces the $(n+i)$th row $w_i$ by $w_i+a v_j$ and the $(n+j)$th row $w_j$ by $w_j-a v_i$. 

Similarly,
\begin{align*}
    G_{ij}(a):= \begin{pmatrix}
                I_n & a(E_{ij}-E_{ji})\\
                0 & I_n
                \end{pmatrix}.
\end{align*}
Multiplying $M$ by $G_{ij}(a)$ from the left replaces the $i$th row $v_i$ by $v_i+a w_j$ and the $j$th row $v_j$ by $v_j-a w_i$.


\subsection{Odd case}
For a matrix $M \in M_{2n+1}(R)$ denote the rows of $M$ by $$v_1,\dots,v_n, v_0, w_1,\dots,w_n. $$
We consider the following \emph{odd elementary orthogonal matrices} in $SO_{2n+1}(R, S_{2n+1})$. As in the even case we have the following three basic types. Let $a \in R$.

For $i\neq j$, set 
\begin{center}
$    K_{ij}(a):=
            \begin{pmatrix}
                I_n+aE_{ij} & 0 & 0\\
                0 & 1 & 0\\
                0 & 0 & I_n-aE_{ji}
            \end{pmatrix}.$
\end{center}

For $1\le i<j\le n$, set
\begin{center}
$F_{ij}(a):=
\begin{pmatrix}
I_n & 0 & 0\\
0 & 1 & 0\\
a(E_{ij}-E_{ji}) & 0 & I_n
\end{pmatrix}, \quad
G_{ij}(a):=
\begin{pmatrix}
I_n & 0 & a(E_{ij}-E_{ji}) \\
0 & 1 & 0\\
0 & 0 & I_n
\end{pmatrix}.$
\end{center}

In addition, there are two new types of elementary orthogonal matrices involving the middle row and the middle column:
\begin{center}
 $   P_i(a)=
    \begin{pmatrix}
        I_n&ae_i&-\dfrac{a^2}{2} E_{ii}\\
        0&1&-ae_i^T\\
        0&0&I_n
    \end{pmatrix} \quad \text{and} \quad
    Q_i(a)=
    \begin{pmatrix}
        I_n   &   0 &   0\\
        a e_i^T   &   1 &   0\\
        -\dfrac{a^2}{2} E_{ii}  &   -a e_i    &   I_n
    \end{pmatrix}.$
\end{center}
Multiplying $M$ by $P_i(a)$ from the left replaces the row $v_j$ by $v_j + a \delta_{ij} v_0 - \frac{a^2}{2} \delta_{ij} w_i$, the row $v_0$ by $v_0 - a w_i$ and leaves the row $w_j$ unchanged. 
Multiplying $M$ by $Q_i(a)$ from the left leaves the row $v_j$ unchanged, replaces the row $v_0$ by $v_0 + a v_i$ and the row $w_j$ by $w_j - a \delta_{ij} v_0 - \frac{a^2}{2} \delta_{ij} v_i$. 

\begin{remark}
    To put things in the frame of unitriangular factorizations, elementary factors of type $F$ and $G$ are lower unitriangular and upper unitriangular respectively. Factors of type $K_{ij}$ with $i < j$ becomes upper unitriangular under the change of basis vectors
    \[
        b_1, \dots, b_n, b_{n+1}, \dots, b_{2n} \to  b_1, \dots, b_n, b_{2n}, \dots, b_{n+1}. 
    \]
    Compare with the discussion in \cite[pp.\@ 2-3]{MR5101717}.
\end{remark}




\subsection{Whitehead lemma}\label{subsect-Whitehead}
In the Gauss--Jordan elimination and the Gram--Schmidt process, we need the multiplication with the following diagonal matrices. Consider for $a \in R^*$ 
\begin{center}
$    K_{ii}(a)=\begin{pmatrix}
                    \mathrm{diag}(1,\dots,a,\dots,1) & 0\\
                    0 & \mathrm{diag}(1,\dots,a^{-1},\dots,1)
                \end{pmatrix} $  in the even case 
\end{center}
and 
\begin{center}
  $  K_{ii}(a)=
\begin{pmatrix}
\mathrm{diag}(1,\dots,a,\dots,1) & 0 & 0\\
0 & 1 & 0\\
0 & 0 & \mathrm{diag}(1,\dots,a^{-1},\dots,1)
\end{pmatrix}$  in the odd case.
\end{center}

These matrices are however not unipotent. For our unitriangular factorization it is thus critical to know if $K_{ii}(a)$ is a product of elementary orthogonal matrices.  Such a factorization of a diagonal matrix was called the Whitehead lemma by Bass \cite{MR0249491} and studied for the symplectic group in \cite{MR4078081}. Here we will see $K_{ii}(a)$ can be written as a product of elementary orthogonal factors if $a$ is a square in $R$. Namely there is $b \in R$ such that $a = b^2$. In particular over the complex numbers $K_{ii}(a)$ can be factorized for all $a \in \C^*$. 

We may restrict to the even case since the computations leave the middle row and the middle column intact for the odd case. 
The following formula, well-known to K-theory specialists, appeared e.g.\@ in Vavilov \cite[Section 5]{05859446} in the form of semisimple long roots of Lie algebra. Here we present it as  multiplication of group elements. 


\begin{proposition}[Whitehead lemma]
\label{prop:whitehead} 
Fix $1\le i\le n$. Let $a\in R^\ast$ be a square with $t^2=a$ for some $t\in R^\ast$.
Then for every $1 \le j \neq i \le n$, 
\begin{align*}
 K_{ii}(a) =  w^-_{ij}(t)\,w^-_{ij}(1)^{-1} w^+_{ij}(t)\,w^+_{ij}(1)^{-1} 
\end{align*}
where 
\begin{center}
   $ w^-_{ij}(t):=K_{ij}(t)\,K_{ji}(-t^{-1})\,K_{ij}(t) \,$
    and $\,
    w^+_{ij}(t):=G_{ij}(t)\,F_{ij}(t^{-1})\,G_{ij}(t).$
\end{center}
In particular $K_{ii}(a)$
is a product of elementary orthogonal matrices of the type $K_{ij}$, $F_{ij}$, $G_{ij}$.
\end{proposition}
\begin{proof}
    Direct computation shows
    \begin{align*}
    w^-_{ij}(t)
    &=
    K_{ij}(t)\,K_{ji}(-t^{-1})\,K_{ij}(t) \\
    &=
    \begin{pmatrix}
    I_n-E_{ii}-E_{jj}-t^{-1}E_{ji}+tE_{ij} & 0\\
    0 & I_n-E_{ii}-E_{jj}-tE_{ji}+t^{-1}E_{ij}
    \end{pmatrix}
    \end{align*}
    and 
    \begin{align*}
        w^+_{ij}(t) &= G_{ij}(t)\,F_{ij}(t^{-1})\,G_{ij}(t) 
    =
    \begin{pmatrix}
    I_n-E_{ii}-E_{jj} & t(E_{ij}-E_{ji})\\
    t^{-1}(E_{ij}-E_{ji}) & I_n-E_{ii}-E_{jj}
    \end{pmatrix}.
    \end{align*}
    Then
    \begin{center}
    $w^-_{ij}(t) w^-_{ij}(1)^{-1}=K_{ii}(t)\,K_{jj}(t^{-1}) \quad \text{ and } \quad w^+_{ij}(t) w^+_{ij}(1)^{-1} =K_{ii}(t)\,K_{jj}(t).$
    \end{center}
    Their product is $K_{ii}(t^2) = K_{ii}(a)$ since diagonal matrices commute. 
\end{proof}

\begin{remark} \label{rema-whitehead-odd}
   (a) For $SO_m(R,S_m)$ with $m \geq 6$ and $R$ commutative with $2 \in R^{\ast}$, the elementary factors of type $K$ can be completely replaced by type $F$ and $G$. For $K_{ij}$ with $i \neq j$ we have 
   $$K_{ij}(a) = G_{ik}(a)F_{kj}(1)G_{ik}(-a)F_{kj}(-1), \quad k \neq i, j.$$

    (b) Combining (a) with Proposition \ref{prop:whitehead}, we essentially only need elementary orthogonal matrices of type $F$ and $G$ for the factorization when $m \ge 6$ is even. In such cases we get an orthogonal counterpart of the factorization for $Sp_{2n}(\C)$ as in  \cite{MR4078081}, where only elementary symplectic types similar to $F$ and $G$ are allowed. 
    
    For $m= 4$ a factorization of $K_{ij}$ as in (a) is not possible, in contrast to the symplectic case. This is due to the fact that $\mathfrak{so}_4 = \mathfrak{sl}_2 \oplus \mathfrak{sl}_2$, where the $K$-type roots span one $\mathfrak{sl}_2$ and the $F$ and $G$ type span the other. In the symplectic case the Lie algebra $\mathfrak{sp}_{2n}$ is simple. 

   (c) In the odd case, $K_{ii}(a)$ can also be factorized into a product of only $P_i$ and $Q_i$. Fix $1\le i\le n$ and let $a\in R^\ast$ be a square with $t^2=a$ for some $t\in R^\ast$. Let
    \begin{center}
    $w_{i}(t):=P_{i}(t)\,Q_{i}(-2t^{-1})\,P_{i}(t).$
    \end{center}
    We get
    \begin{center}
    $K_{ii}(a)=w_{i}(t)\,w_{i}(1)^{-1}.$
    \end{center}
    This also covers $m=3$, where Proposition \ref{prop:whitehead} does not apply since there is no $j\neq i \in \{ 1 \}$.

    (d) Moreover, for odd $m \ge 5$ and $i \neq j$, since 
    \begin{align*}
        [P_i(a), Q_j(b)] = K_{ij}(ab), \quad [P_i(a), P_j(b)] = G_{ij}(-ab), \quad [Q_i(a), Q_j(b)] = F_{ij}(-ab),
    \end{align*}
    a factorization in the style of \cite{MR4078081} with only $P_i, Q_i$ as elementary orthogonal matrices is possible. 

    For $m=3$, the only elementary
orthogonal matrices are $P_1$ and $Q_1$. We already have from (c) that $K_{11}(a)$ is a product of them.
\end{remark}

\section{The Gauss--Jordan Elimination near the identity} \label{sec-GJ}
Let $R$ be a commutative unital Banach algebra. On its matrix algebra $M_m(R)$ one has the norm
$\|M\|=\max_{p,q}\|M_{pq}\|$, and on $SO_m(R, S_m)$ the induced topology. 
\begin{proposition}
\label{lemm-GJ-orth}
Assume $m\geq 3$. There exists a neighborhood $V$ of the identity in $SO_{m}(R,S_m)$ and an integer $N=N(m)$ such that
every $M\in V$ can be written as
\begin{center}
$M = E_1(M)\cdots E_N(M),$
\end{center}
where each $E_\nu(M)$ is of the form $K_{ii}$ or one of the elementary matrices
$K_{ij},F_{ij},G_{ij}$ (and $P_i, Q_i$ for the odd case) with parameters depending continuously on $M\in V$,
and $E_\nu(I_m)=I_m$ for all $\nu$.
\end{proposition}

\begin{proof}
\textbf{Even case}.
Let 
$$
    M=\begin{pmatrix}A&B\\ C&D\end{pmatrix} \in SO_{2n}(R,S_{2n}).
$$
If $M$ is sufficiently close to $I_{2n}$, then $A$  is close to $I_n$ and each pivot during the elimination stays invertible after shrinking $V$, since it is continuous in $M$ and equals $1$ at $I_{2n}$.
Then, we can reduce $A$ to $I_n$ by using left multiplication by $K_{ij}$ and $K_{ii}$.
Here, all parameters depend continuously on $M$, as they are built from matrix entries of $M$ and their inverses, which are close to $1$.

Now, we multiply from the left by $F_{ij}$ to make the lower-left block $C$ lower triangular.
Since $A=I_n$, we have that the first $n$ entries of the row $v_k$ is the $k$-th standard basis vector of $\mathbb{C}^n$.
Hence, the row operations coming from multiplying $F_{ij}$ allow us to eliminate off-diagonal entries of $C$
as in the Gauss--Jordan process. After finitely many steps, we make $C$ lower triangular.

For $M\in SO_{2n}(R,S_{2n})$, the relation $M^TS_{2n}M=S_{2n}$ gives the identity $A^T C + C^T A = 0$.
Because $A=I_n$, we have $C+C^T=0$. Since a skew-symmetric lower triangular matrix
must be zero, we get $C = 0$. 

The remaining block identities become
\begin{center}
$A^T D + C^T B = I_n,\ B^T D + D^T B = 0.$
\end{center}
With $A=I_n$ and $C=0$ we get $D=I_n$ and $B^T+B=0$, so $B$ is skew-symmetric.
Finally, with $C=0$ and $D=I_n$, multiplying by $G_{ij}$ adds a skew elementary matrix to the block $B$
while keeping $C$ and $D$ fixed. Hence we can make $B$ upper triangular by Gaussian elimination using
$G_{ij}$. Since $B$ remains skew-symmetric, this forces $B=0$.

Thus $M$ is reduced to $I_{2n}$ by a fixed finite sequence of elementary steps, and reversing the steps yields the
claimed factorization on some neighborhood $V$ of $I_{2n}$. 

\textbf{Odd case}. 
We reduce it to the even case. Write
\begin{center}
$M=
\begin{pmatrix}
A&u&B\\
x&\alpha&y\\
C&v&D
\end{pmatrix}, 
    \,\, \text{ where } \alpha\in R,\
    x,y\in R^{1\times n},\
    u,v\in R^{n\times 1},
    A,B,C,D\in M_n(R).$
\end{center}

The identity $M^T S_{2n+1} M= S_{2n+1}$ gives
\begin{align}
\alpha^2+2u^Tv&=1,                                   \label{equa-odd-gj-1}
\\
\alpha x+u^TC+v^TA&=0,                                  \label{equa-odd-gj-2}\\
x^Tx+A^TC+C^TA&=0,                                      \label{equa-odd-gj-3}\\
x^Ty+A^TD+C^TB&=I_n,                                    \label{equa-odd-gj-4}\\
\alpha y+u^TD+v^TB&=0,                                  \label{equa-odd-gj-5}\\
y^Ty+B^TD+D^TB&=0.    \nonumber                                  \label{equa-odd-gj-6}
\end{align}
After every elementary operation we retain the same notation for
the blocks of the resulting matrix.

As in the proof for the even case, we first reduce $A$ to $I_n$, using left multiplication by
$K_{ij}$ and $K_{ii}$ to perform 
Gauss--Jordan elimination on $A$.
More explicitly, at the $k$-th stage let $\lambda_k$ be the 
$(k,k)$-entry of $A$. First multiply $K_{kk}(\lambda_k^{-1})$ from the left
to normalize the  $(k,k)$-entry. Then for every $i\neq k$, multiply $K_{ik}(-a_{ik})$ to eliminate other entries in the $k$-th column of $A$. 

After this elimination process, we have
\begin{center}
$M=
\begin{pmatrix}
I_n&u&B\\
x&\alpha&y\\
C&v&D
\end{pmatrix}.$
\end{center}

We then use $P_{i}$ and $Q_{i}$ to eliminate $x,y,u,v$. 

For $a\in R$, recall that
\begin{center}
$P_i(a)=
\begin{pmatrix}
I_n&ae_i&-\dfrac{a^2}{2}e_ie_i^T\\
0&1&-ae_i^T\\
0&0&I_n
\end{pmatrix}.$
\end{center}
Then
\begin{equation*}
MP_i(a)=
\begin{pmatrix}
I_n
&
u+ae_i
&
B-au e_i^T-\dfrac{a^2}{2}e_ie_i^T
\\
x
&
\alpha+ax_i
&
y-\left(a\alpha+\dfrac{a^2}{2}x_i\right)e_i^T
\\
C
&
v+aCe_i
&
D-av e_i^T-\dfrac{a^2}{2}Ce_ie_i^T
\end{pmatrix}.
\end{equation*}
Choosing $a = - u_i$ and multiplying $M$ with $P_i$ for $i = 1, 2, \dots, n$ we get $u=0$. 

The matrix now has the form
\begin{center}
    $M=
\begin{pmatrix}
I_n&0&B\\
x&\alpha&y\\
C&v&D
\end{pmatrix}.$
\end{center}
Applying left multiplication by $Q_i(-x_i)$, $i = 1, 2, \dots, n$ successively we get $x= 0$. 

With $\ u=0,\ A=I_n,\ x=0,$ Equation \eqref{equa-odd-gj-2} gives $v=0$. Since $\det M = 1 \neq 0$, the middle column is nonzero and thus $\alpha \neq 0$. Moreover, $\alpha^2=1$ by Equation \eqref{equa-odd-gj-1} and so $\alpha$ is a unit. 
By Equation \eqref{equa-odd-gj-5}, we also get $y=0$. The matrix $M$ is now of the form 
\begin{center}
$M=
\begin{pmatrix}
I_n&0&B\\
0&\alpha&0\\
C&0&D
\end{pmatrix}.$
\end{center}
This implies
\begin{center}
$\det M = \alpha \det \begin{pmatrix} I_n &B\\ C&D\end{pmatrix}$.
\end{center}
By \eqref{equa-odd-gj-3} and \eqref{equa-odd-gj-4} we have $D-CB = I_n$. From Schur's complement formula, 
$$
    \det \begin{pmatrix} I_n &B\\ C&D\end{pmatrix} = \det( D - CB) = 1, 
$$
thus $\alpha = \det M=1$. The middle column and the middle row are now decoupled from the rest and we may proceed as before in the even case.
\end{proof}

\begin{remark} \label{rmk-K_ii}
    After the Gauss--Jordan elimination near the identity, some diagonal factors $K_{ii}(a)$ close to the identity remain. 
    By the Whitehead lemma, Proposition \ref{prop:whitehead} (Remark \ref{rema-whitehead-odd}(c) when $m=3$), these matrices can be factorized into products of elementary matrices. Indeed, since we have $||a - 1|| < \epsilon$ with $\epsilon < 1$,
    \[
    \sqrt{a} := \exp (\frac{1}{2} \log a ), \quad \log a:= - \sum_{j \ge 1} \frac{(1-a)^j}{j}.
    \]
\end{remark}

\section{nullhomotopy factorization over Banach algebras} \label{sec-Banach}

In this section $R$ is a commutative Banach algebra with unit.
We write
\begin{center}
$SO_{m}(R,S_m):=\{M\in SL_{m}(R): M^T S_m M=S_m \}.$
\end{center}
We equip $SO_{m}(R,S_m)$ with the metric induced by a matrix norm on $M_{m}(R)$.

\begin{theorem}
\label{thm-main}
Let $m \geq 3$ and $M\in SO_{m}(R,S_m)$ be nullhomotopic
in $SO_{m}(R,S_m)$. Then $M$ is a finite product of elementary orthogonal matrices of the form
$K_{ij},F_{ij},G_{ij}$ (and $P_i, Q_i$ for the odd case) with parameters in $R$.
\end{theorem}
\begin{proof}
Let $M_t$ be a nullhomotopy in $SO_{m}(R,S_m)$ with $M_0=I_m$ and $M_1=M$.
By uniform continuity of $t\mapsto M_t$ and continuity of inversion, there exists $\delta>0$ such that
$M_tM_{t'}^{-1}\in V$ whenever $|t-t'|<\delta$, where $V$ is as in Proposition~\ref{lemm-GJ-orth}.
Choose an integer $k>1/\delta$. Then
\begin{center}
$M=M_1=\bigl(M_1M_{1-\frac1k}^{-1}\bigr)\bigl(M_{1-\frac1k}M_{1-\frac2k}^{-1}\bigr)\cdots
\bigl(M_{\frac1k}M_0^{-1}\bigr),$
\end{center}
a product of $k$ elements of $V$. Each factor is a product of elementary matrices by Proposition~\ref{lemm-GJ-orth} and Remark \ref{rmk-K_ii},
hence so is $M$.
\end{proof}

As a corollary, we have 
\begin{corollary}

Let $m \ge 3$. For a compact  Hausdorff topological space $X$, every nullhomotopic continuous map
$X\to SO_{m}(\mathbb{C}, S_m)$  admits
a factorization into elementary orthogonal factors with entries in $\mathcal{C}(X)$.
\end{corollary}

\section{The Gram--Schmidt process and the uniform nullhomotopy} \label{sec-GS}
Let $m \ge 3$. We use the standard Hermitian inner product on $\mathbb C^{m}$,
\begin{center}
$\langle x,y\rangle := x\,\overline{y}^{\,T},\ \|x\|^2 := \langle x,x\rangle.$
\end{center}

\begin{proposition}
\label{theo-GS}
There is $L = L(m)$ such that for every $M\in SO_{m}(\mathbb{C}, S_m)$ there are orthogonal matrices $E_1, \dots,  E_L \in SO_{m}(\mathbb{C}, S_m)$  of the form $K_{ij}$, $K_{ii}$, $F_{ij}$ or $Q_i$ 
such that
\begin{center}
$(E_1\cdots E_L)M \in SO_{m}(\mathbb{C}, S_m)\cap U(m).$
\end{center}
\end{proposition}

\noindent \textit{Proof of Proposition} \ref{theo-GS}.
For the even case $m = 2n$ denote by $v_1,\dots,v_n,w_1,\dots,w_n$ the rows of $M$. 
We proceed in three steps, keeping the notation $v_1,\dots,v_n,w_1,\dots,w_n$ for the rows of the current matrix.

\paragraph{\textbf{Step 1}: we make $v_1,\dots,v_n$ pairwise orthogonal.}
For $k=1,2,\dots,n-1$ and $i=k+1,\dots,n$ multiply $M$ from the left by
\begin{center}
$K_{ik}\!\left(-\frac{\langle v_i,v_k\rangle}{\|v_k\|^2}\right).$
\end{center}
This gives the replacement
\begin{center}
$v_i\leftarrow v_i-\frac{\langle v_i,v_k\rangle}{\|v_k\|^2}\,v_k,$
\end{center}
which removes the component of $v_i$ along $v_k$. Proceeding inductively makes $v_1,\dots,v_n$
pairwise orthogonal.

\paragraph{\textbf{Step 2}: we normalize $v_1,\dots,v_n$.}
Left multiplications by $K_{ii}(1/\|v_i\|),\ i=1,\dots,n,$
make $\|v_i\|=1$ for all $i$. 

\paragraph{\textbf{Step 3}: we make the first $n$ rows orthogonal to the last $n$ rows.}
Let
\begin{center}
$B=\bigl(\langle v_i,w_j\rangle\bigr)_{1\le i,j\le n}.$
\end{center}
For each fixed $j=1,\dots,n-1$ and each $i=j+1,\dots,n$, multiply by
\begin{center}
    $F_{ji}\!\left(-\langle w_j,v_i\rangle\right).$
\end{center}
This replaces $w_j\leftarrow w_j-\langle w_j,v_i\rangle v_i$ and the new $w_j$ is orthogonal to $v_i$. After finitely many steps, the matrix $B$ becomes upper triangular.

After Steps~1-3, we claim the resulting matrix lies in $U(2n)$.
Consider the Hermitian Gram matrix $H:=MM^\ast.$
Then one can show that $H\in SO_{2n}(\mathbb{C}, S_{2n})$ by using $M \in SO_{2n}(\mathbb{C}, S_{2n})$, $S_{2n}^{-1} = S_{2n}$ and $\overline{S}_{2n} = S_{2n}$. 



Write $H$ in $n\times n$ blocks:
\begin{center}
$H=\begin{pmatrix}A&B\\ C&D\end{pmatrix}.$
\end{center}
By the previous steps, the $v$'s are orthonormal, hence $A=I_n$. By Step~3, $B$ is upper triangular.
Since $H$ is Hermitian, $C=B^\ast$ is lower triangular. Finally, the orthogonal identity $H^TSH=S$
implies the block relations
\begin{align}
A^T C + C^T A &= 0, \label{equa-H1}\\
A^T D + C^T B &= I_n. \label{equa-H2}
\end{align}
Since $A=I_n$, \eqref{equa-H1} becomes $C+C^T=0$, i.e.\ $C$ is skew-symmetric. But $C$ is lower triangular,
hence $C=0$. Consequently $B=C^\ast=0$ and \eqref{equa-H2} gives $D=I_n$. Therefore $H= MM^\ast=I_{2n}$,
i.e.\  $M\in U(2n)$.

\medskip
For the odd case $m = 2n+1$ denote  by $v_1, \dots, v_n, v_0, w_1, \dots, w_n$ the row vectors of $M$.

As for the even case, we first use elementary matrices $K_{ij}$ and $K_{ii}$ to achieve that
$\langle v_i,v_j\rangle=\delta_{ij}$ for
$ 1\leq i,j\leq n$. Then we use elementary matrices of type $Q$ to make $v_0$ orthogonal to $v_1,
\dots,v_n$:

\begin{center}
    $\langle v_0,v_i\rangle=0\ \text{for every }1\leq i\leq n.$
\end{center}
The Gram matrix $H =MM^*$ has the form
\begin{center}
    $H=
\begin{pmatrix}
I_n&0&B\\
0&a&c^*\\
B^* &c&D
\end{pmatrix},$
\end{center}
where
$
a=\|v_0\|^2>0$. Solving for the $(2,1)$- and $(2,2)$-block of the orthogonal condition $H^T S_{2n+1} H = S_{2n+1}$ we have $a = 1$ and 
\begin{center}
$c^T =
\begin{pmatrix}
\langle w_1,v_0\rangle &
\cdots &
\langle w_n,v_0\rangle
\end{pmatrix} = 0 .$
\end{center}
That is $v_0$ is already orthogonal to all the $w_j$'s. Hence, $H$ simplifies to
\begin{center}
    $H=
\begin{pmatrix}
I_n&0&B\\
0&1&0\\
B^*&0&D
\end{pmatrix}.$
\end{center}
Now the Gram--Schmidt process for the even case applies and   the proof is complete. \hfill \qed

\begin{remark} \label{cont-GS}
    For the ring of complex-valued continuous functions on a normal topological space $X$, the above Gram--Schmidt process also works. Given a continuous $M : X \to SO_m(\C, S_m)$, the coefficients in the above elementary orthogonal matrices $E_j(M)$ depend continuously on $M$. At $M = I_m$ the rows are orthonormal and so every parameter is $0$ or $1$, namely $E_j(I_m) = I_m$. 
\end{remark}

Recall
\begin{center}
    $SO_{m}(\mathbb C):=\{A\in SL_{m}(\mathbb C): A^T A=I_{m}\}.$
\end{center}

\begin{lemma} \label{lem:model}
One has $SO_{m}(\mathbb C)\cap U(m)=SO_{m}(\mathbb R)$.
\end{lemma}

\begin{proof}
Let $A\in SO_{m}(\mathbb C)\cap U(m)$. Then $A^TA=I_m$ and $A^\ast A=I_m$, hence $A^T=A^{-1}=A^\ast$.
Comparing entries of $A^T$ and $A^\ast$ gives $A_{ji}=\overline{A}_{ji}$ for all $i,j$, so $A$ has real entries,
i.e.\ $A\in SO_{m}(\mathbb R)$. The reverse inclusion is clear.
\end{proof}

Define
\begin{center}
    $P_0 := \frac{1}{\sqrt2}\begin{pmatrix} I_n & iI_n\\ I_n & -iI_n\end{pmatrix} \text{ for } m = 2n$ and $P_0:=\frac{1}{\sqrt2}\begin{pmatrix} I_n & 0 & iI_n\\ 0 & \sqrt2 & 0\\ I_n & 0 & -iI_n\end{pmatrix}$ for $m= 2n+1$.
\end{center}
\begin{lemma}
The matrix $P_0$ lies in $U(m)$ and satisfies $P_0^T S_m P_0 = I_{m}$ and $\det(P_0)=(-i)^n$.
\end{lemma}

\begin{proof}
A direct computation.
\end{proof}

\begin{lemma} \label{lemma-conj}
For $M\in SL_{m}(\mathbb C)$ set $A := P_0^{-1} M P_0$. Then $M^T S_m M = S_m$ if and only if $A^T A = I_{m}$. 
In particular, conjugation by $P_0$ gives an isomorphism of algebraic groups
\begin{center}
$SO_{m}(\mathbb{C}, S_m)\xrightarrow{\ \sim\ } SO_{m}(\mathbb C),\
M\longmapsto P_0^{-1}MP_0 .$
\end{center}
\end{lemma}

\begin{proof}
Assume $M^T S_m M=S_m$ and write $M=P_0 A P_0^{-1}$. Then
\begin{center}
$S_m = M^T S_m M = (P_0^{-T} A^T P_0^T)\,S_m \,(P_0 A P_0^{-1})
= P_0^{-T}\,A^T\,(P_0^T S_m P_0)\,A\,P_0^{-1}.$
\end{center}
Multiplying by $P_0^T$ on the left and by $P_0$ on the right gives
\begin{center}
$P_0^T S_m P_0 = A^T (P_0^T S_m P_0) A.$
\end{center}
Because $P_0^T S_m P_0=I_{m}$, we then get $A^T A=I_{m}$. For the converse implication, we obtain it by reversing the computation.
\end{proof}

\begin{proposition} \label{K_U-compact}
Conjugation by $P_0$ induces an isomorphism of topological groups
\begin{center}
$K_U = SO_{m}(\mathbb{C}, S_m)\cap U(m)\xrightarrow{\ \sim\ } SO_{m}(\mathbb R),\
M\longmapsto P_0^{-1}MP_0 .$
\end{center}
In particular, $K_U$ is a compact subgroup of $SO_{m}(\mathbb{C}, S_m)$.
\end{proposition}

\begin{proof}
Let $M\in K_U$ and let $A=P_0^{-1}MP_0$. Since $P_0\in U(m)$, the conjugation preserves $U(m)$. Hence,
$A\in U(m)$. Also, $M\in SO_{m}(\mathbb{C}, S)$ implies $A\in SO_{m}(\mathbb C)$ by the previous lemma.
Thus, $A\in SO_{m}(\mathbb C)\cap U(m)=SO_{m}(\mathbb R)$.

Conversely, if $A\in SO_{m}(\mathbb R)$ then $A\in U(m)$ and $A^T A=I_m$, so $M=P_0 A P_0^{-1}$ lies in
$U(m)$ and in $SO_{m}(\mathbb{C}, S_m)$. Hence $M \in K_U$ and the map is bijective with continuous inverse,
so it is a homeomorphism. Therefore $K_U$ is compact.
\end{proof}
Since $SO_{m}(\mathbb{R}), m \ge 3$ is compact with finite fundamental group $\mathbb{Z}/2\mathbb{Z}$, we obtain the following by applying \cite[Theorem 1]{MR0458416} , cf.\ \cite{03707400} and \cite[Proof of Theorem 2]{04073958}.

\begin{proposition}\label{Cardel-Siegel}
Let X be a finite dimensional normal space and
assume $M : X \longrightarrow K_{U}$ is nullhomotopic. Then there exists a uniform homotopy $M_t : X \longrightarrow K_{U}$ with $M_1 = M$ and $M_0 \equiv I_m$, i.e.\ for any $\varepsilon > 0$ there is $\delta > 0$ such that $||M_t(x) - M_{t^{\prime}}(x)|| < \varepsilon$ for all $x \in X$ whenever $|t - t^{\prime}| < \delta$.
\end{proposition}

\section{Proof of main theorem} \label{sec:combine}

\begin{proof}[Proof of Theorem \ref{cont-factorization}]

First, for $m=3$, we use Remark \ref{rema-whitehead-odd}(c) instead of Proposition \ref{prop:whitehead}, otherwise the proof is the same for all $m \ge 3$.

Let $M_t$ be the nullhomotopy with $M_t: X \to SO_{m} (\mathbb C, S_m)$ with $M_1 = M$ and $M_0 = I_m$. By the Gram--Schmidt process in Proposition \ref{theo-GS} and Remark \ref{cont-GS}, there are elementary orthogonal matrices
$E_1, \cdots, E_L$ such that 
$$
    M^{\prime}_t = E_1(M_t)E_2(M_t)\cdots E_L(M_t)M_t
$$ 
is a nullhomotopy of continuous maps from $X$ to $SO_m(\C, S_m)$ with values in $K_U$. By Proposition \ref{Cardel-Siegel} there is a uniform nullhomotopy $\tilde{M}_t: X \to K_U$ with
$\tilde{M}_1 = E_1(M)\cdots E_L(M)M$. Let $V$ be the neighborhood of the identity in $SO_{m} (\mathbb C, S_m)$ from Proposition \ref{lemm-GJ-orth} with $R = \mathbb C$. Since $K_U$ is contained in $U(m)$, there is $\varepsilon > 0$ such that $AB^{-1} \in V$ for all $A, B \in K_U$ with $\|A - B\| < \varepsilon$. Hence, for $\delta$ as in Proposition \ref{Cardel-Siegel} and an integer $k > 1/\delta$, the maps of form 
\begin{center}
    $D_j(x) := \tilde{M}_{j/k}(x)\, \tilde{M}_{(j-1)/k}(x)^{-1}, \ j = 1, \dots, k,$
\end{center}
take values in $V$, and $\tilde{M}_1(x) = D_k(x) \cdots D_1(x)$. Applying Proposition \ref{lemm-GJ-orth} to each $D_j$, 
we get $\tilde{M}_1(x) = E^{\prime}_1(x)\cdots E^{\prime}_l(x)$ for all $x \in X$. With $l = k N(m)$  all values of $D_j(x)$ lie in the neighborhood $V$ as in Proposition \ref{lemm-GJ-orth}. Thus we get a factorization
\begin{center}
 $M(x) = E^{-1}_L(x) \cdots E^{-1}_1(x)E^{\prime}_1(x)\cdots E^{\prime}_l(x)$
\end{center}
with all the elementary factors depending continuously on $x \in X$. Now replace factors of the form $K_{ii}$ by elementary factors: first the $K_{ii}$s from the Gauss--Jordan elimination factorize by Remark \ref{rmk-K_ii}, whereas the $K_{ii}$ from the normalization step in the Gram--Schmidt process have positive parameter and thus Proposition \ref{prop:whitehead} applies. 

Finally, the existence of a uniform bound on the number of factors depending only on $m$ and $d$ follows Vaserstein's argument in \cite[proof of Theorem 2]{04073958}. Suppose such a uniform bound does not exist, then we could choose normal topological spaces $X_i$ of covering dimension at most $d$ and nullhomotopic map $f_i : X_i \to SO_m(\C, S_m)$ that need more than $i$ factors. Then one takes the disjoint union $X = \sqcup_i X_i$ which is normal of covering dimension at most $d$. The map $f$ with $f \vert_{X_i} = f_i$ is nullhomotopic, thus $f$ is a product of some number $N$ factors. Restricting $f$ to $X_{N+1}$ gives a contradiction. 
\end{proof}

\section*{Acknowledgement}
We would like to thank Frank Kutzschebauch for helpful discussions and comments on an earlier draft.

\end{document}